\documentclass[11pt]{pjm}

\def\R{{\hbox{\bf R}}}
\def\FRAME{{\mathcal F}}

\def\p{{\hbox{\bf p}}}

\font \roman = cmr10 at 10 true pt

\def\Im{{\hbox{\roman Im}}}

\def\be#1{\begin{equation}\label{#1}}
\def\bas{\begin{align*}}
\def\eas{\end{align*}}
\def\bi{\begin{itemize}}
\def\ei{\end{itemize}}

\def\dim{{\hbox{\roman dim}}}

\def\Z{{\hbox{\bf Z}}}
\def\eps{\varepsilon}

\def\emph#1{{\it #1}}
\def\textbf#1{{\bf #1}}

\def\gain{{\frac{1}{16}}}

\theoremstyle{plain}
  \newtheorem{theorem}[subsection]{Theorem}

  \newtheorem{proposition}[subsection]{Proposition}
  
  \newtheorem{lemma}[subsection]{Lemma}
  \newtheorem{corollary}[subsection]{Corollary}

\theoremstyle{remark}

\theoremstyle{definition}
  \newtheorem{definition}[subsection]{Definition}

\include{psfig}

\begin{document}

\title[Besicovitch sets in four dimensions]{A new bound for finite field Besicovitch sets in four dimensions}

\author{Terence Tao}
\address{Department of Mathematics, UCLA, Los Angeles CA 90095-1555}
\email{tao@@math.ucla.edu}
\urladdr{http://www.math.ucla.edu/~tao}
\thanks{ This work was inspired by the April 2002 Instructional Conference on Combinatorial aspects of Mathematical Analysis at University of Edinburgh.  The author is a Clay Prize Fellow and is supported by the Packard Foundation.  The author thanks Tony Carbery, Nets Katz, Wilhelm Schlag, and Jim Wright for helpful discussions, and is also indebted to David Gieseker and Allen Knutson for their explanation of some of the basics of algebraic geometry.  The author is particularly indebted to Nets Katz for emphasizing the importance of reguli to this problem.  Finally, the author thanks the anonymous referee for careful reading and detection of several misprints.} 


\begin{abstract}  Let $F$ be a finite field with characteristic greater than two.  Define a \emph{Besicovitch set} in $F^4$ to be a set $P \subseteq F^4$ containing a line in every direction.  The \emph{Kakeya conjecture} asserts that $|P| \approx |F|^4$.  In \cite{wolff:survey} it was shown that $|P| \gtrsim |F|^3$.  In this paper we improve this to $|P| \gtrapprox |F|^{3+\gain}$. On the other hand, we show that the bound of $|F|^3$ is sharp if we relax the assumption that the lines point in different directions.  One new feature in the argument is the introduction of a small amount of basic algebraic geometry.
\end{abstract}

\maketitle

\section{Introduction}

Let $F$ be a finite field with characteristic greater than 2.  For any $n \geq 2$, we define a \emph{Besicovitch set} in $F^n$ to be a set $P \subseteq F^n$ containing a line in every direction.  The \emph{finite field Kakeya conjecture} (see e.g. \cite{wolff:survey}) asserts that $|P| \geq C_\eps |F|^{n-\eps}$ for any $\eps > 0$, where $|P|$ denotes the cardinality of $P$ and the quantities $C_\eps$ are independent of $|F|$.
This conjecture is the finite field analogue of the Euclidean Kakeya set conjecture, which is related to several other problems in harmonic analysis; see 
\cite{wolff:survey}, \cite{gerd:kakeya} for further discussion on this.  Basically, one can view the finite field Kakeya problem as a simplified model problem for the more interesting Euclidean Kakeya problem, where several technical difficulties (involving small separations, small angles, or multiple scales) have been eliminated.

Informally, the Kakeya conjecture asserts that lines which point in different
directions in $F^n$ cannot have substantial overlap.
This conjecture has been proven in two dimensions but remains open in higher dimensions.  In \cite{wolff:survey} (see also \cite{wolff:kakeya}, \cite{gerd:kakeya}) it was shown that $|P| \gtrsim |F|^{(n+2)/2}$ (here $A \gtrsim B$ denotes the estimate $A \geq C^{-1} B$).  In fact, more was proven:

\begin{definition}\label{wolff-def}  A family $L$ of lines in $P^n$ is said to obey the \emph{Wolff axiom} if for every $2 \leq k \leq n-1$, every $k$-dimensional affine subspace\footnote{By affine subspace we mean a translation of a vector subspace of $F^n$.} $V \subset F^n$ contains at most $O(|F|^{k-1})$ lines in $L$.
\end{definition}

\begin{theorem}\label{wolff-bound}\cite{wolff:kakeya}, \cite{wolff:survey}  If $L$ is a family of $O(|F|^{n-1})$ lines obeying the Wolff axiom, and $P \subseteq F^n$ contains all the lines in $L$, then $|P| \gtrsim |F|^{(n+2)/2}$.
\end{theorem}

In fact one only needs to use the Wolff axiom for $k=2$.  From this theorem and the observation that any family of lines which point in different directions
automatically obeys the Wolff axiom, we immediately see that Besicovitch sets have cardinality $\gtrsim |F|^{(n+2)/2}$.

In \cite{gerd:kakeya} (see also \cite{katzlaba}) it was observed that the statement of Theorem \ref{wolff-bound} was sharp in three dimensions, in the sense that there existed finite fields $F$ and collections of lines $L$ in $F^3$ obeying the Wolff axiom and a collection $P$ of points containing all the lines in $L$, such that $|P| \sim |F|^{5/2}$.  Indeed, if $F$ contained a subfield $G$ of index 2, with the accompanying involution $z \mapsto \overline z$ on $F$, then one could take $P$ to be the \emph{Heisenberg group}
$$ P := \{ (z_1, z_2, z_3) \in F^3: \Im(z_3) = \Im(z_1 \overline z_2) \}$$
where $\Im(z) := (z - \overline z)/2$.  (It is an interesting question as to
whether an example similar to this could be obtained if $F$ did not
contain a subfield of index 2).

Our first observation is that Theorem \ref{wolff-bound} is also sharp in four-dimensions:

\begin{proposition}\label{prop}  Let $\langle , \rangle: F^4 \times F^4 \to F$ be a non-degenerate symmetric quadratic form on $F^4$.  Let $P$ be the ``unit sphere''
\be{sphere}
P := \{ x \in F^4: \langle x, x\rangle = 1 \}
\end{equation}
and let $L$ be the set of all lines of the form $\{ x+tv: t \in F \}$, where $x \in F^4$, $v \in F^4 \backslash \{0\}$ are such that $\langle x, x\rangle = 1$, $\langle v, x \rangle = 0$, and $\langle v, v\rangle = 0$.  Then we have that $L$ has cardinality $|L| \sim |F|^3$ and obeys the Wolff axiom, while $P$ has cardinality $|P| \sim |F|^3$ and contains all the lines in $L$.
\end{proposition}

We prove this in Section \ref{counter-sec}. A similar counterexample can be created in $\R^4$ as long as one chooses the signature of the form $\langle,\rangle$ to be indefinite.  Observe that the above Proposition does not contradict the Kakeya conjecture because the lines $L$ do not all point in different directions (despite obeying the Wolff axiom).  Nevertheless, it seems of interest to extend this example (and the Heisenberg group) to higher dimensions, though perhaps the bound of $|F|^{(n+2)/2}$ in Theorem \ref{wolff-bound} is not necessarily sharp for large $n$.

This example illustrates two things.  Firstly, in order to progress toward the Kakeya conjecture in low dimensions\footnote{In high dimensions (e.g. $n \geq 9$) there are other, more ``arithmetic'' arguments available to improve upon Theorem \ref{wolff-bound}. See \cite{borg:high-dim}, \cite{KT}, \cite{kt:kakeya2}, \cite{rogers}, \cite{gerd:kakeya}.} one must make more use the hypothesis that the lines in $L$ point in different directions; merely assuming the Wolff axiom will not by itself suffice.  Secondly, the algebraic geometry of quadric surfaces may be relevant to the Kakeya problem\footnote{There seems to be a parallel phenomenon in recent work on Szemer\'edi's theorem on arithmetic progressions, in that while arithmetic progressions are rather ``linear'' quantities, they give rise rather naturally to other ``quadratic'' objects which then need to be studied.  See \cite{gowers}.}.

In the three-dimensional case $n=3$, quadric surfaces are essentially the same thing as \emph{reguli} - the ruled surfaces consisting of all the lines which intersect three fixed lines in general position.  In particular, we have the ``three-line lemma'', which asserts that given three mutually skew lines in $F^3$, there are at most $O(|F|)$ lines in different directions which intersect all three.

Reguli have already come up in the work of Schlag \cite{schlag:kakeya}, who used the three-line lemma to give a new proof of Bourgain's estimate \cite{borg:kakeya}
\be{b-est}
|P| \gtrsim |F|^{7/3}
\end{equation} in three dimensions.  While it is true that this bound has since been superceded by Wolff's estimate in Theorem \ref{wolff-bound}, we shall need to follow \cite{schlag:kakeya} and make use of reguli and the three-line lemma in what follows.  We are indebted to Nets Katz for pointing out the usefulness of reguli in the low-dimensional Kakeya problem.  Indeed, our work here was inspired by similar work in three dimensions by Nets Katz (currently in preparation).

The main result of this paper is the following improved bound on the cardinality of Besicovitch sets in four dimensions.  We use $A \lessapprox B$ to denote the estimate $A \leq C_\eps |F|^\eps B$ for any $\eps > 0$, where $C_\eps$ is a quantity depending only on $\eps$.

\begin{theorem}\label{main}  If $P$ is a Besicovitch set in $|F|^4$, then $|P| \gtrapprox |F|^{3 + \gain}$.
\end{theorem}

One can probably improve the $\gtrapprox$ to a $\gtrsim$ by going through
the argument in this paper more carefully, but we will not 
do so here in order to simplify the exposition.

The paper is organized as follows.  After setting out our incidence geometry notation in Section \ref{incidence-sec}, we prove Proposition \ref{prop} in Section \ref{counter-sec}.  We then review some basic algebraic geometry in Section \ref{algebraic-sec}, culminating in a ``three-regulus lemma'' in $F^4$, which will be the analogue of the three-line lemma in $F^3$.  We then review some combinatorial preliminaries in Section \ref{notation-sec}, before starting the proof of Theorem \ref{main} in Section \ref{start-sec}.  The first step is to use a standard ``iterated popularity'' argument (as in e.g. \cite{ccc}), together with a rudimentary version of the ``plate number'' argument in \cite{wolff:xray}, in order to refine the Besicovitch set to a uniform, non-degenerate collection of points and lines.  Once we have done a sufficient number of refinements, we can construct a large number of reguli which are incident to many lines in the Besicovitch set, and eventually get about $|F|^3$ lines incident to three distinct reguli (if $|P|$ is too close to $|F|^3$); this will contradict the three regulus 
lemma mentioned earlier.

\section{Incidence notation}\label{incidence-sec}

We now set some notation for the finite field geometry of the affine space $F^4$.  A \emph{line} in $F^4$ is a set of the form $l = \{ x + tv: t \in F \}$ where $x, v \in F^4$ and $v$ is non-zero.  Two lines are \emph{parallel} if they are translates of each other but not identical; a set of lines is said to \emph{point in different directions} if no two lines in the set are parallel or identical. 

A \emph{2-plane} in $F^4$ is a set of the form $\pi = \{ x + t_1 v_1 + t_2 v_2: t_1, t_2 \in F \}$ where $x,v_1,v_2 \in F^4$ and $v_1$, $v_2$ are linearly independent.  Two lines are \emph{coplanar} if they lie in the same 2-plane; observe that coplanar lines must either be identical, parallel, or intersect in a point.  A pair of lines are \emph{skew} if they are not coplanar.

A \emph{3-space} in $F^4$ is a set of the form $\lambda = \{ x + t_1 v_1 + t_2 v_2 + t_3 v_3: t_1,t_2,t_3 \in F \}$ where $x,v_1,v_2,v_3 \in F^4$ and $v_1$, $v_2$, $v_3$ are linearly independent.  Observe that any pair of skew lines lies in a unique 3-space.  Two 3-spaces are \emph{parallel} if they are disjoint, and one is the
translate of the other.

We shall use the symbol $p$ to refer to points, $l$ to lines, $\pi$ to 2-planes, and $\lambda$ to 3-spaces.  We use the symbol $P$ to refer to sets of points, $L$ to sets of lines, $\Pi$ to sets of 2-planes, and $\Lambda$ to sets of 3-spaces.  We use $Gr(F^4, 1)$ to denote the space of all lines, $Gr(F^4, 2)$ to denote the space of all 2-planes, and $Gr(F^4, 3)$ to denote the space of all 3-spaces.  (Note that these are the \emph{affine} Grassmanians, in that the spaces do not need to contain the origin).

\section{The counterexample}\label{counter-sec}

We now prove Proposition \ref{prop}.  It is likely that this Proposition follows from the standard theory of Fano varieties of quadric surfaces, but we will just give an elementary argument.

Let $P$ and $L$ be as in Proposition \ref{prop}.  It is clear from construction that the lines in $L$ lie in $P$.  Now we verify the cardinality bounds.
We begin with a standard lemma on the number of ways of representing a field
element as a quadratic form.

\begin{lemma}\label{sum}  Let $\langle, \rangle: F^n \times F^n \to F$ be a symmetric bilinear form on $F^n$ with rank at least 1, and let $Q(x) := \langle x, x \rangle$ be the associated quadratic form.  Then we have
\be{ss}
\{ (x_1, \ldots, x_n) \in F^n: Q(x_1, \ldots, x_n) = x \} \lesssim |F|^{n-1}
\end{equation}
for all $x \in F$.  

If we know that $\langle,\rangle$ has rank at least 3, then we can improve \eqref{ss} to
\be{ss-2}
\{ (x_1, \ldots, x_n) \in F^n: Q(x_1, \ldots, x_n) \rangle = x \} \sim |F|^{n-1}
\end{equation}
for all $x \in F$, if $|F|$ is sufficiently large.  
\end{lemma}

\begin{proof}
By placing the quadratic form $Q$ in normal form (recalling that ${\rm char} F \neq 2$) we may assume that
$$ Q(x_1, \ldots, x_n) = \alpha_1 x_1^2 + \ldots + \alpha_k x_k^2$$
where $k$ is the rank of $Q$ and $\alpha_1, \ldots, \alpha_k$ are non-zero elements of $F$.  We may assume that $k=n$ since the general case $n \geq k$ follows by adding $n-k$ dummy variables.  In particular $\alpha_j \neq 0$ for $j=1, \ldots, n$.

The bound \eqref{ss} is now clear, since if we fix $x_1, \ldots, x_{n-1}$ and $x$ then there are at most 2 choices for $x_n$.  Now let us assume $n \geq 3$, and prove \eqref{ss-2}.

We use Gauss sums. We fix $e: F \to S^1$ to be a non-principal character of $F$, i.e. a multiplicative function from $F$ to the unit circle which is not identically 1. For instance, if $F = \Z/p\Z$ for some prime $p$, one can take 
$e(x) := \exp(2\pi i x/p)$.

For any $y \in F$, let $S(y)$ be the Gauss sum
$ S(y) := \sum_{x \in F} e(y x^2)$.
As is well known (see e.g. \cite{gerd:kakeya}), $S(0)$ is equal to $|F|$, while $|S(y)| = |F|^{1/2}$ for all non-zero values of $y$.  

Fix $x \in F$.  By expanding the Kronecker delta as a Fourier series,
we see that the number of solutions to \eqref{ss} can be written as
\bas 
\sum_{x_1, \ldots, x_n \in F} \delta(\alpha_1 x_1^2 + \ldots + \alpha_n x_n^2 - x)
&= \frac{1}{|F|} \sum_{y \in F} 
\sum_{x_1, \ldots, x_n \in F} e((\alpha_1 x_1^2 + \ldots + \alpha_n x_n^2 - x)y)\\
&= \frac{1}{|F|} \sum_{y \in F} 
e(-xy) \prod_{i=1}^n S(\alpha_i y)\\
&= |F|^{n-1} + \frac{1}{|F|} \sum_{y \in F \backslash \{0\}} 
e(-xy) \prod_{i=1}^n S(\alpha_i y)\\
&=  |F|^{n-1} + \frac{1}{|F|} \sum_{y \in F \backslash \{0\}} 
O(|F|^{n/2})\\
&= |F|^{n-1} + O(|F|^{n/2})
\end{align*}
as desired, since $n \geq 3$.
\end{proof}

From the Lemma we see that $|P| \sim |F|^3$, as desired.  Now we count the lines in $L$.  

From the Lemma we see that there are $\sim |F|^3$ choices of null direction $\{v \in F^4 \backslash 0: \langle v, v\rangle = 0 \}$.  For each such $v$, the space $v^\perp := \{ x \in F^4: \langle x, v \rangle = 0 \}$ is 3-dimensional (since $Q$ is non-degenerate).  Furthermore, since $Q$ is non-degenerate on $F^4$ and $v$ is a null vector, we see that $Q$ must also be non-degenerate on $v^\perp$.  Restricting $Q$ to $v^\perp$ (which is of course isomorphic to $F^3$) we see from Lemma \ref{sum} that there are $\sim |F|^2$ choices for $x$.  Thus there are $\sim |F|^5$ possible pairs $(x,v)$ that generate a line in $L$.  However, each line in $L$ is generated by $\sim |F|^2$ such pairs $(x,v)$, so we have $|L| \sim |F|^3$ as desired.

It remains to verify the Wolff axiom.  First pick a 3-space $\lambda$ and consider the lines in $L$ which go through $\lambda$.  

Pick an arbitrary point $x_0$ in $\lambda$, so that $\lambda-x_0$ is a three-dimensional subspace of $F^4$.  By Lemma \ref{sum}, the number of null vectors $\{ v \in \lambda - x_0: \langle v, v \rangle = 0 \}$ is $O(|F|^2)$.  

Fix $v$ as above.  There are two cases.  If $v^\perp \not \equiv (\lambda-x_0)$, then there are $O(|F|)$ choices of $x \in \lambda$ such that $\langle x, v \rangle = 0$ and $\langle x, x \rangle = 1$.  But if $v^\perp \equiv (\lambda - x_0)$, then the number of choices for $x$ could be as large as $O(|F|^2)$. But $(\lambda - x_0)^\perp$ only has cardinality $O(|F|)$, hence the number of $v$ in the second category is at most $O(|F|)$.  Thus the number of pairs $(x,v)$ which can generate a line in $\lambda$ is at most $O(|F|^3)$. But each line is generated by $\sim |F|^2$ pairs $(x,v)$.  Thus the number of lines in $\lambda$ is at most $O(|F|)$, which clearly implies the Wolff axiom for both $k=2$ and $k=3$.  This completes the proof of Proposition \ref{prop}.

\section{Some basic algebraic geometry}\label{algebraic-sec}

Here we review some basic facts from algebraic geometry (see e.g. \cite{harris}), and apply them to our Kakeya problem.  The material we will need is not very advanced; basically, we need the concept of the dimension of an algebraic variety, and we need to know that this dimension behaves in the expected way with respect to intersections, projections, cardinality, etc.  We shall also rely heavily on the basic fact that the dimension of an algebraic variety is always an \emph{integer} (in contrast to, say, the ``half-dimensional'' field $G$ mentioned in the introduction).

Let $\overline{F}$ denote the algebraic closure of $F$ and $n \geq 1$.  An \emph{algebraic variety} in $\overline{F}^n$ is defined to be the zero locus of a collection $Q_1, \ldots, Q_k$ of $\overline{F}$-valued polynomials on the affine space $\overline{F}^n$.  In this paper we shall always assume that our algebraic varieties have bounded degree, thus $k=O(1)$ and all the polynomials $Q_1, \ldots, Q_k$ have degree $O(1)$.  

An algebraic variety $V$ in $\overline{F}^n$ has a well-defined \emph{integer-valued} dimension $0 \leq d \leq n$; there are several equivalent definitions of this dimension, for instance $d$ is the smallest non-negative integer such that generic affine spaces in $\overline{F}^n$ of codimension greater than $d$ are disjoint from $V$.  (See \cite{harris} for more equivalent definitions of dimension).  
If $V$ has dimension $n$ then it must be all of $\overline{F}^n$, while if $V$ has dimension 0 then it can only consist of at most $O(1)$ points. Of course, the algebraic geometry notion of dimension is consistent with the linear algebra notion of dimension, thus for instance 3-spaces have dimension 3.

An algebraic variety is \emph{irreducible} if it does not contain any proper sub-variety of the same dimension.  Every algebraic variety of dimension $k$ can be decomposed as a union of $O(1)$ irreducible varieties of dimension at most $k$ (see e.g. \cite{harris}).

We define an \emph{algebraic variety} in $F^n$ of dimension $d$ to be a restriction to $F^n$ of an algebraic variety in $\overline{F}^n$ of dimension $d$.  Observe that if $V$ is a variety in $F^n$ of dimension $d$ then $|V| \lesssim |F|^d$ (this can be shown, for instance, by taking generic intersections with affine spaces of codimension $d$).

Let $L \subseteq Gr(F^n, 1)$ be a collection of lines which point in different directions.  In the introduction we observed that this implies the Wolff axiom, that not too many lines in $L$ can lie inside a $k$-space.  In fact we can generalize this to $k$-dimensional varieties (cf. \cite{gerd:kakeya}, Proposition 8.1):

\begin{lemma}[Generalized Wolff property]\label{wolff-variety} Let $V \subseteq F^n$ be an algebraic variety in $F^n$ of dimension $k$, and let $L \subseteq Gr(F^n, 1)$ be a collection of lines in $F^n$ which point in different directions.  Then we have
$$ | \{ l \in L: l \subseteq V \} | \lesssim |F|^{k-1}.$$
\end{lemma}

Observe that the lines in Proposition \ref{prop} violate the above property, but of course those lines do not point in different directions\footnote{On the other hand, one can show that the lines arising from the Heisenberg example do obey this generalized Wolff property.  It may be that in the three-dimensional theory, one needs to extend this lemma further, to cover not only varieties over $F$, but also over subfields of $F$ such as $G$.}.

\begin{proof}  
We may of course assume that $|F| \gg 1$, since the claim is obvious for
$|F|$ bounded.

We can embed $F^n$ in the projective space $PF^{n+1}$, which we think of as the union of $F^n$ with the hyperplane at infinity.  By replacing the defining polynomials of $V$ with their homogeneous counterparts, we can thus extend $V$ to a $k$-dimensional variety $\overline{V}$ in $PF^{n+1}$ (see e.g. \cite{harris}).  

We break up $\overline{V}$ into irreducible components, each of dimension at most $k$.  We can assume that none of the irreducible components are contained inside the hyperplane at infinity, since we could simply remove those components and still have an extension of $V$.  In particular we see that the intersection of $\overline{V}$ with the hyperplane at infinity is at most $k-1$-dimensional.

Let $l$ be a line in $L$, which we can extend to be a projective line $\overline{l}$ in $PF^{n+1}$ by adding a single point at infinity (the direction of $l$).  Observe that the restriction of $\overline{V}$ to $\overline{l}$ is an algebraic variety of dimension either 0 or 1; in other words, either the projective line $\overline{l}$ lies inside $\overline{V}$, or else $\overline{l}$ intersects $\overline{V}$ in at most $O(1)$ points.  Thus in order for $l$ to be contained in $V$, the direction of $l$ must lie inside $\overline{V}$ (assuming that $|F|$ is sufficiently large).  But by the previous paragraph the number of such directions is at most $O(|F|^{k-1})$.  Since the lines in $L$ point in different directions, we are done.
\end{proof}

As a consequence of this lemma we see that a Besicovitch set cannot have high
intersection with algebraic variety:

\begin{corollary}\label{variety}  Let $V \subseteq F^n$ be an algebraic variety in $F^n$ of dimension at most $n-1$, and let $L \subseteq Gr(F^n, 1)$ be a collection of lines in $F^n$ which point in different directions.  Then we have
$$ | \{ (p,l) \in V \times L: p \in l \}| \lesssim |F|^{n-1}.$$
\end{corollary}

Note that the trivial upper bound for the left-hand side is $|F| |L| \lesssim |F|^n$.  Thus this lemma gains a power of $|F|$ over the trivial bound.

\begin{proof} As in the proof of Lemma \ref{wolff-variety}, we observe that every line $l$ in $F^n$ is either contained in $V$, or else intersects $V$ in at most $O(1)$ points. The lines of the second type contribute at most $O(|L|) = O(|F|^{n-1})$ incidences, while by Lemma \ref{wolff-variety} the lines of the first type contribute at most $O(|F| |F|^{\dim(V)-1}) = O(|F|^{n-1})$ incidences, and we are done.
\end{proof}

A further consequence is that the \emph{lines} of a Besicovitch set cannot have
large intersection with an algebraic variety:

\begin{corollary}\label{lines} Let $L \subseteq Gr(F^n, 1)$ be a collection of lines in $F^n$ which point in different directions, and let $P \subseteq F^n$ be a set of points containing all the lines in $L$.  Let $W \subseteq Gr(F^n,1)$ be an algebraic variety \emph{of lines} of dimension at most $n-1$.  Then we have
$$ | L \cap W | \lesssim |F|^{n-2} + |F|^{-1} |P|.$$
\end{corollary}

Again, this lemma gains a power of $|F|$ over the trivial bound of $|F|^{n-1}$ (assuming $P$ is not too huge).

\begin{proof}  Consider the set
$$ X := \{ (p,l) \in F^n \times W: p \in l \}.$$
This is an algebraic variety in $F^n \times Gr(F^n,1)$ of dimension at most $n$.  Now consider the map $\phi: X \to F^n$ given by $\phi(p,l) := p$.  Observe that for any $p$ in the image of $\phi$, the fibers $\phi^{-1}(p)$ are either 0-dimensional (i.e. have cardinality $O(1)$), or at least 1-dimensional.  This implies
(see e.g. \cite{harris}) that we have a decomposition $\phi(X) := P_1 \cup P_2$, where the fibers $\phi^{-1}(p)$ are 0-dimensional for all $p \in P_1$, and $P_2$ is contained in an algebraic variety of dimension at most $n-1$.

By construction of $P_1$ we have
\bas | \{ (p, l) \in P_1 \times (L \cap W): p \in l \} |
&\lesssim | \{ p \in P_1: p \in l \hbox{ for some } l \in L \cap W \} |\\
&\lesssim | \{ p \in P_1: p \in P \} |\\
& \leq |P|.
\end{align*}
Also, by Corollary \ref{variety} we have
$$ | \{ (p, l) \in P_2 \times (L \cap W): p \in l \} | \lesssim |F|^{n-1}.$$
Adding the two estimates, we see that
$$ |F| |L \cap W| = | \{ (p,l) \in \phi(X) \times (L \cap W): p \in l \}| \lesssim |F|^{n-1} + |P|$$
and the claim follows.
\end{proof}

To apply the above results to our four-dimensional
problem, we need some notation for reguli.

\begin{definition}\label{frame-def}  A \emph{frame} $f$ is a quadruplet $f = (l_1, l_2, l_3, \lambda)$ where $\lambda \in Gr(F^4, 3)$ is a 3-space in $F^4$, and $l_1, l_2, l_3 \in Gr(F^4, 1)$ are distinct, mutually skew lines in $F^4$ which lie inside $\lambda$.  If $f$ is a frame, we write $\lambda(f)$ for $\lambda$.  If $f = (l_1, l_2, l_3, \lambda)$ is a frame, we use $L(f)$ to denote the set of lines $l \in Gr(F^4, 1)$ which intersect $l_1$, $l_2$, and $l_3$, and $r(f) \subseteq \lambda$ to denote the union of all the lines in $L(f)$.
\end{definition}

We refer to $r(f)$ as the \emph{regulus} generated by a frame $f$.  It is well known (see e.g. \cite{harris}) that a regulus is a quadric surface\footnote{A model example of a regulus is the set $\{ (x,y,xy,0) \in F^4: x,y \in F \}$, where the lines $l_i$ are of the form $\{ (x, y_i, xy_i, 0): x \in F \}$ for some distinct $y_1, y_2, y_3 \in F$.  In fact all reguli can be shown to be projectively equivalent to this example.} in $\lambda$, i.e. it is the zero locus of some quadratic polynomial in $\lambda$.  In particular it is an algebraic variety of dimension 2.  Since the lines in a frame are mutually skew, we see that this quadratic polynomial is irreducible (so the regulus is not a (double) plane, or the union of two planes), and that the lines $L(f)$ have cardinality $\sim |F|$ and are finitely overlapping. 

\begin{corollary}[Three regulus lemma]\label{3-reg} Let $L \subseteq Gr(F^4, 1)$ be a collection of lines in $F^4$ which point in different directions, and let $P \subseteq F^4$ be a set of points containing all the lines in $L$.  Let $f_1, f_2, f_3$ be three frames such that the 3-spaces $\lambda(f_1)$, $\lambda(f_2)$, $\lambda(f_3)$ are parallel and disjoint.  Then
$$ | \{ l \in L: l \cap r(f_i) \neq \emptyset \hbox{ for all } i = 1,2,3 \} |
\lesssim |F|^2 + |F|^{-1} |P|.$$
\end{corollary}

Again, note that this bound improves by roughly $|F|$ over the trivial bound of $|F|^3$, if $|P|$ is not much larger than $|F|^3$.  The hypothesis that the 3-spaces $\lambda(f_i)$ are parallel can be substantially relaxed, but we will not need to do so here.

\begin{proof}
Fix $f_1, f_2, f_3$, and let $W \subseteq Gr(F^4, 1)$ denote the set
\be{W-def}
W := \{ l \in Gr(F^4, 1): l \cap r(f_i) \neq \emptyset \hbox{ for all } i = 1,2,3 \}.
\end{equation}
Since the $r(f_i)$ and $l$ are algebraic varieties, it is clear (e.g. by using resultants; see e.g. \cite{harris}) that the relationship $l \cap r(f_i) \neq \emptyset$ is equivalent to some finite set of explicit algebraic relations between the defining parameters of $l$ and $f_i$.  Thus $W$ is an algebraic variety in $Gr(F^4, 1)$.  In light of Corollary \ref{lines}, it will suffice to verify that $W$ has dimension\footnote{We apologize to readers expert in algebraic geometry for the appalling crudeness of the following argument.} at most 3.

Let $p$ be a point in $r(f_1)$.  Let $\phi_p$ be the stereographic projection from $\lambda(f_2)$ to $\lambda(f_3)$, thus $\phi_p(x) = y$ iff $p$, $x$, $y$ are collinear.  Observe that $W$ is isomorphic to
$$ \{ (p, y) \in r(f_1) \times \lambda(f_3): y \in r(f_3) \cap \phi_p(r(f_2)) \}$$
(basically because two points determine a line, and because the planes $\lambda(f_i)$ are disjoint).  In other words, one can think of $W$ as a bundle over $r(f_1)$ whose fiber at $p$ is $r(f_3) \cap \phi_p(r(f_2))$.

Note that $\phi_p: \lambda(f_2) \to \lambda(f_3)$ is an invertible linear map, so that $\phi_p(r(f_2))$ is an irreducible quadric surface in $\lambda(f_3)$.  The set $r(f_3) \cap \phi_p(r(f_2))$ thus has dimension at most 2, and in fact will have dimension at most 1 unless $\phi_p(r(f_2)) \equiv r(f_3)$ (by irreducibility).  However, as $p$ varies, the quadric surfaces $\phi_p(r(f_2))$ move by translation.  Since $r(f_2)$ is not a plane, we thus see that there can be at most a one-dimensional family of points $p$ for which $\phi_p(r(f_2)) \equiv r(f_3)$.  

To summarize, as $p$ varies over the two-dimensional variety $r(f_1)$, the fiber $r(f_3) \cap \phi_p(r(f_2))$ is at most one-dimensional, except possibly for a one-dimensional family of points $p$ for which the fiber is two-dimensional.  From this it is clear\footnote{Strictly speaking, to compute the dimension properly one should work in the algebraically closed field $\overline{F}$ here.  But this causes no difficulty as the above geometric considerations are valid for all fields of characteristic larger than two.} that $W$ has dimension at most 3, and we are done.
\end{proof}

Corollary \ref{3-reg} is the analogue of the three lines lemma used in \cite{schlag:kakeya}.  Our strategy will now be to start with a Besicovitch set and construct many frames $f$ and many lines $l \in L$ so that $r(f)$ intersects $L$, in order to exploit the above Corollary.  To do this we shall need some basic combinatorial tools, which we now pause to review.

\section{Some basic combinatorics}\label{notation-sec}

We shall frequently use the following elementary observation:  If $B$ is a finite set, and $\mu: B \to \R^+$ is a function such that
$$ \sum_{b \in B} \mu(b) \geq X,$$
then we have
$$ \sum_{b \in B: \mu(b) \geq X/2|B|} \mu(b) \geq X/2.$$
We refer to this as a ``popularity'' argument, since we are restricting $B$ to the values $b$ which are ``popular'' in the sense that $\mu$ is large.  We shall in fact iterate this popularity argument a large number of times.

We shall frequently use the following version of the Cauchy-Schwarz and H\"older inequalities.

\begin{lemma}\label{cz}  Let $A$, $B$ be finite sets, 
and let $\sim$ be a relation connecting pairs $(a,b) \in A \times B$ such that
$$ |\{ (a,b) \in A \times B: a \sim b \}| \gtrsim X$$
for some $X \gg |B|$.  Then
$$ |\{ (a,a',b) \in A \times A \times B: a \neq a'; a, a' \sim b \}| \gtrsim \frac{X^2}{|B|}$$
and
$$ |\{ (a,a',a'',b) \in A \times A \times B: a,a',a'' \hbox{ distinct}; a,a',a'' \sim b \}| \gtrsim \frac{X^3}{|B|^2}$$
\end{lemma}

\begin{proof}
Define for each $b \in B$, define $\mu(b) := | \{ a \in A: a \sim b \} |$.  Then by hypothesis we have
$$ \sum_{b \in B} \mu(b) \gtrsim X.$$
In particular, by the popularity argument we have
$$ \sum_{b \in B: \mu(b) \gtrsim X/|B|} \mu(b) \gtrsim X.$$
By hypothesis, we have $X/|B| \gg 1$.  From this and the previous, we obtain
$$ \sum_{b \in B: \mu(b) \gtrsim X/|B|} \mu(b) (\mu(b)-1) \gtrsim X (X/|B|)$$
and
$$ \sum_{b \in B: \mu(b) \gtrsim X/|B|} \mu(b) (\mu(b)-1) (\mu(b)-2) \gtrsim X (X/|B|) (X/|B|).$$
The claims follow. 
\end{proof}

A typical application of the above Lemma is the standard incidence bound

\begin{corollary}\label{easy-cor}  For an arbitrarily collection $P \subseteq F^n$ of points and $L \subseteq Gr(F^n, 1)$ of lines, we have
\be{easy-incidence}
|\{ (p,l) \in P \times L: p \in l \}| \lesssim |P|^{1/2} |L| + |P| 
\end{equation}
\end{corollary}

\begin{proof}
We may of course assume that the left-hand side of \eqref{easy-incidence} is $\gg |P|$, since the claim is trivial otherwise.  From Lemma \ref{cz} we have
$$ |\{ (p, l, l') \in P \times L \times L: p \in l \cap l'; l \neq l' \}| \gtrsim
|P|^{-1} |\{ (p,l) \in P \times L: p \in l \}|^2.$$
On the other hand, $|l \cap l'|$ has cardinality $O(1)$ if $l \neq l'$, thus
$$ |\{ (p, l, l') \in P \times L \times L: p \in l \cap l'; l \neq l' \}| \lesssim 
|L|^2.$$
Combining the two estimates we obtain the result.
\end{proof}

The above estimate will be most useful when $|L|$ is small (in particular if $|L| = O(|F|)$).  When $|L|$ is large, we have the alternate estimate 

\begin{proposition}\label{wip}\cite{gerd:kakeya}  Let the notation be as in the previous Corollary.  If we further assume that the lines in $L$ point in different directions, then we have
\be{wolff-incidence}
|\{ (p,l) \in P \times L: p \in l \}| \lesssim |P|^{1/2} |L|^{3/4} |F|^{1/4} + |P| + |L|.
\end{equation}
\end{proposition}

\begin{proof}
See \cite{gerd:kakeya}, Proposition 8.6; the argument there is essentially due to Nets Katz, but the original result of this type dates back to Wolff \cite{wolff:kakeya}, \cite{wolff:survey}. 

For the convenience of the reader we now sketch an informal ``probabilistic'' derivation of \eqref{wolff-incidence}.  Let $I$ denote the set in \eqref{wolff-incidence}.  We may assume that $|I| \gg |P|, |L|$ since the claim is trivial otherwise.  

Observe that a randomly chosen point $p \in P$ and a randomly chosen line $l \in L$ have a probability $|I|/|P||L|$ of being incident (so that $p \in l$).  Thus, given two random lines $l_1, l_2 \in L$ and a random point $p \in P$, we expect\footnote{This of course assumes independence of various random events, which is usually not the case.  To make the argument rigorous one must use such tools as Lemma \ref{cz}, which can be viewed as a statement that certain events are positively correlated.  See \cite{gerd:kakeya}, Proposition 8.6 for details.} the chance that $p$ is incident to both $l_1$ and $l_2$ is $(|I|/|P||L|)^2$.  Since there are $|P|$ possible values for $p$, the chance that two random lines $l_1, l_2 \in L$ intersect at all is thus heuristically $|P| (|I|/|P||L|)^2$.

As a consequence, the probability that three random lines $l_1, l_2, l_3 \in L$ form a triangle is heuristically $(|P| (|I|/|P||L|)^2)^3$.  (There is the chance that this triangle is degenerate, but the hypothesis $|I| \gg |P|, |L|$ can be used to show that the probability of this occurring is low).  On the other hand, given two intersecting lines $l_1, l_2 \in L$, there are at most $O(|F|)$ lines $l_3 \in L$ which can intersect them both, since we may apply the Wolff axiom to the 2-plane spanned by $l_1$ and $l_2$.  Combining these estimates we obtain
$$ (|P| (|I|/|P||L|)^2)^3 \lesssim (|P| (|I|/|P||L|)^2)  \frac{|F|}{|L|}$$
and \eqref{wolff-incidence} follows\footnote{It is an instructive exercise to obtain similar heuristic probabilistic derivations of such estimates as \eqref{easy-incidence} (using the fact that two random lines intersect in at most one point) or \eqref{b-est} (using the fact that a random regulus contains at most $O(|F|)$ lines).  See also Section \ref{heuristic-sec} below.}.
\end{proof}

We remark that one only requires the Wolff axiom on $L$ to obtain \eqref{wolff-incidence}. In particular one can easily obtain Theorem \ref{wolff-bound} as a consequence of \eqref{wolff-incidence}. It is likely that one can generalize Theorem \ref{main} to obtain a further improvement to \eqref{wolff-incidence}, but we do not pursue this question here. 

\section{A heuristic proof of Theorem \ref{main}}\label{heuristic-sec}

We now give a heuristic explanation as to why we can improve upon Theorem \ref{wolff-bound} in four dimensions, in the spirit of the probabilistic arguments in Proposition \ref{wip}.  In later sections we shall make this heuristic argument rigorous.

Suppose for contradiction we have a family $L \subseteq Gr(F^4, 1)$ of lines in different directions of cardinality $|L| \sim |F|^3$ which are contained in a set $P \subseteq F^4$, also of cardinality $|P| \sim |F|^3$. Arguing as in Proposition \ref{wip} we see that any two lines in $L$ have a (heuristic) probability $\sim 1/|F|$ of intersecting.

Also, a random line $l \in L$ and a random 3-space $\lambda \in Gr(F^4, 3)$ have a probability $1/|F|^2$ of being incident (so that $l \subseteq \lambda$).  Thus we expect a 3-space $\lambda$ to contain $|L|/|F|^2 \sim |F|$ lines in $L$.

Now consider the set of all quintuples $(l_1, l_2, l^1, l^2, l^3) \in L^5$ of lines such that $l_i$ intersects $l^j$ for all $i=1,2$, $j=1,2,3$.  From the above heuristics we see that there should be about $|L|^5 (1/|F|)^6 \sim |F|^9$ such quintuples.  On the other hand, for generic quintuples $(l_1, l_2, l^1, l^2, l^3)$ of the above form, the lines $l^1$, $l^2$, $l^3$ must lie in a 3-space $\lambda$, and $l_1$, $l_2$ must lie in the regulus generated by the frame $(l^1, l^2, l^3, \lambda)$.  (For this heuristic argument we ignore the possibility that the quintuple could degenerate).

To count the number of possible reguli, observe that there are $|L|^2 \sim |F|^6$ choices for $l^1, l^2$, which determines $\lambda$.  From our previous heuristic we see that $\lambda$ can contain at most $O(|F|)$ choices for $l^3$, thus there are at most $O(|F|^7)$ reguli.

Dividing $|F|^9$ by $|F|^7$, we thus see that a generic regulus $r(f)$ of the above type must contain at least $|F|^2$ pairs $(l_1, l_2)$ of lines in $L$.  But $r(f)$ only has $O(|F|)$ lines to begin with.  Thus a generic regulus $r(f)$ must have extremely large intersection with $P$, so that $|r(f) \cap P| \sim |r(f)| \sim |F|^2$.

Since a random $p \in P$ and $l \in L$ have a probability $1/|F|^2$ of being incident, this means that a random line $l \in L$ and a random regulus $r(f)$ have a probability $\sim 1$ of intersecting.  In particular, if we select three parallel reguli $r(f_1)$, $r(f_2)$, $r(f_3)$, a large fraction of lines in $L$ must be incident to all three reguli.  But this contradicts Corollary \ref{3-reg}, since $|L| \sim |F|^3$ and $|P| \sim |F|^3$.

\section{Preliminary refinements}\label{start-sec}

We now begin the rigorous proof of Theorem \ref{main}, which will broadly follow the heuristic outline of the previous section.

Let $P_0 \subseteq F^4$ be a Besicovitch set.  We may assume that
\be{p0-small}
|P_0| \lessapprox |F|^{3 + \gain}.
\end{equation}
since the claim is trivial otherwise.  We may also assume that $|F| \gg 1$ for similar reasons. 
 
Since $P_0$ is a Besicovitch set, there exists a set $L_0 \subseteq Gr(F^4,1)$ of lines in different directions such that $|L_0| \sim |F|^3$ and $P_0$ contains every line in $L_0$.  In particular the incidence set
$$ I_0 := \{ (p,l) \in P_0 \times L_0: p \in l \}$$
has cardinality $|I_0| = |F| |L_0| \sim |F|^4$.

Given any line $l$ in $L_0$ and a randomly selected 3-space $\lambda$ in $Gr(F^4,3)$, the probability that $l$ lies in $\lambda$ is $\sim 1/|F|^2$.  Since $|L_0| \sim |F|^3$, one thus expects every 3-space $\lambda$ contains about $|F|$ lines in $L_0$ on the average.  A similar heuristic leads us to expect every 2-plane $\pi \in Gr(F^4,2)$ to contain at most $O(1)$ lines on the average. 

Although these statements need not be true for all 3-spaces $\lambda$, certain variants do hold if we refine $L_0$ and $P_0$ slightly:

\begin{proposition}\label{p-ref}  There exists a quantity
\be{a-small}
1 \lesssim \alpha \lessapprox N^{\gain}
\end{equation}
and a subset\footnote{In the course of this argument we shall need to refine the set $P_0$ to a slightly smaller set $P_1$, and then further to $P_2$, and similarly refine $L_0$ to $L_1$ and then $L_2$, while also refining some auxiliary sets $H_0$ to $H_1$, and $\FRAME_0$ to $\FRAME_1$ to $\FRAME_2$ to $\FRAME_3$.  These refinements are largely technical and as a first approximation one can view these sets as being essentially the same (although the sets $\FRAME_2$, $\FRAME_3$ are significantly smaller than $\FRAME_0$, $\FRAME_1$).} $P_1$ of $P_0$ and a subset $L_1$ of $L_0$ such that the following properties hold.

\begin{itemize}
\item  (Many incidences)  We have the incidence bound
\be{big-incidence}
| \{ (p,l) \in P_1 \times L_1: p \in l \} | \gtrapprox |F|^4.
\end{equation}
\item (Cardinality and multiplicity bounds)
We have the cardinality bound
\be{p1-bound}
|P_1| \lessapprox \alpha |F|^3
\end{equation}
and the multiplicity bound
\be{high-mult}
| \{ l \in L_1: p \in l \} | \approx \alpha^{-1} |F| 
\end{equation}
for all $p \in P_1$.
\item  (No 3-space degeneracy)  For any 3-space $\lambda \in Gr(F^4, 3)$, we have
\be{no-3space}
| \{ l \in L_1: l \subset \lambda \}| \lessapprox \alpha^2 |F|.
\end{equation}
\item  (No 2-plane degeneracy)  For any 2-plane $\pi \in Gr(F^4, 2)$, we have
\be{no-2plane}
| \{ l \in L_1: l \subset \pi \}| \lessapprox \alpha^4.
\end{equation}
\end{itemize}
\end{proposition}

The quantity $\alpha$ measures the improvement over Wolff's bound $|P_0| \gtrsim |F|^3$.  As one can see from \eqref{a-small}, it is rather close to 1. 

\begin{proof}  We follow standard ``iterated refinement'' arguments (see \cite{wolff:xray}, \cite{laba:xray}, \cite{ccc}, \cite{wright:radon}; our argument here is particularly close to that in \cite{laba:xray}).  The purpose of the iteration is mainly to obtain the property \eqref{no-2plane}.

Define the multiplicity function $\mu_0$ on $P_0$ by
$$ \mu_0(p) := |\{ l \in L_0: p \in l \}|.$$
Then we have
$$ \sum_{p \in P_0} \mu_0(p) = |I_0|.$$
If we divide $\mu_0(p)$ into dyadic ``pigeonholes'' and apply the dyadic pigeonhole principle (observing that $\log |F| \approx 1$), we conclude that there exists a multiplicity $\alpha^{-1} |F| $ such that
$$ \sum_{p \in P_0: \mu_0(p) \sim \alpha^{-1} |F| } \mu_0(p) \approx |I_0| \approx |F|^4.$$
Fix this $\alpha$, and define
$$ P'_0 := \{ p \in P_0: \mu_0(p) \sim \alpha^{-1} |F| \} $$
and
$$ I'_0 := \{ (p,l) \in P'_0 \times L_0: p \in l \} \subseteq I_0.$$
Then by construction, $|I'_0| \gtrapprox |I_0| \sim |F|^4$, and
$$ |P'_0| \approx |I_0| / (\alpha^{-1} |F|) \approx \alpha |F|^3.$$
By \eqref{p0-small} we thus have $\alpha \lessapprox N^{\gain}$.  To get the other half of \eqref{a-small}, we observe from Proposition \ref{wip} that
$$
|I'_0| \lesssim |P'_0|^{1/2} |L_0|^{3/4} |F|^{1/4} + |P'_0| + |L_0|;$$
applying the above estimates, we thus obtain $\alpha \gtrsim 1$.  Thus \eqref{a-small} holds.

Set $N := \log \log |F|$; the point of this choice of $N$ is that both $|F|^{C/N}$ and $C^N$ are $\approx 1$ for any fixed choice of constant $C$.
We shall inductively construct sets
\be{p-seq}
P'_0 =: P^{(0)} \supset P^{(1)} \supset \ldots \supset P^{(N)}
\end{equation}
and
\be{l-seq}
L_0 =: L^{(0)} \supset L^{(1)} \supset \ldots \supset L^{(N)}
\end{equation}
as follows\footnote{The use of such a large number of refinements is of course overkill (one could probably get away with $N = 5$, in fact), but reducing the number of refinements used does not alter the exponent $\gain$, since $F^{C/N}$ and $C^N$ were $\approx 1$ anyway.}.

As indicated above, we set $P^{(0)} := P'_0$ and $L^{(0)} := L_0$.  Now suppose inductively that $P^{(k)}$ and $L^{(k)}$ have already been constructed for some $0 \leq k < N$.  We define the incidence set
$$ I^{(k)} := \{ (p,l) \in P^{(k)} \times L^{(k)}: p \in l \}.$$
Clearly we have
$$ \sum_{l \in L^{(k)}} |l \cap P^{(k)}| = |I^{(k)}|.$$
Thus if we set
$$ L^{(k+1)} := \{ l \in L^{(k)}: |l \cap P^{(k)}| \geq \frac{|I^{(k)}|}{2|L^{(k)}|} \}$$
then by the popularity argument we have
$$ \sum_{l \in L^{(k+1)}} |l \cap P^{(k)}| \geq |I^{(k)}|/2.$$
We rewrite this as
$$ \sum_{p \in P^{(k)}} | \{ l \in L^{(k+1)}: p \in l \} | \geq |I^{(k)}|/2.$$
Thus if we set
$$ P^{(k+1)} := \{ p \in P^{(k)}: | \{ l \in L^{(k+1)}: p \in l \} | \geq \frac{|I^{(k)}|}{4 |P^{(k)}|} \}$$
then by the popularity argument again 
$$ \sum_{p \in P^{(k+1)}} | \{ l \in L^{(k+1)}: p \in l \} | \geq |I^{(k)}|/4$$
or in other words
$$ |I^{(k+1)}| \geq |I^{(k)}|/4.$$
We repeat this construction for $k=0,1,\ldots,N-1$, creating a nested sequence of sets of points \eqref{p-seq} and sets of lines \eqref{l-seq}.  By construction and the fact that $4^{-N} \approx 1$, we clearly have
$$ |I^{(k)}| \approx |I^{(0)}| = |I'_0| \gtrapprox |F|^4$$
for all $k$.  Furthermore, we have 
$$ |P^{(k)}| \leq |P'_0| \lessapprox \alpha |F|^3$$
and
$$ |L^{(k)}| \leq |L_0| \lesssim |F|^3.$$
Thus if we set $P_1 := P^{(N)}$ and $L_1 := L^{(N-1)}$, then it is clear that \eqref{big-incidence}, \eqref{p1-bound}, and \eqref{high-mult} hold.  (To get the upper bound in \eqref{high-mult}, simply bound the left-hand side by $\mu_0(p)$).

It remains only to verify the non-degeneracy conditions \eqref{no-3space}, \eqref{no-2plane}.

We first verify \eqref{no-3space}.  Let $\lambda$ be a 3-space.  Since $\lambda$ is clearly an algebraic variety of dimension 3, we can invoke Corollary \ref{variety} and conclude that
$$ |\{ (p,l) \in \lambda \times L^{(N-2)}: p \in l \}| \lesssim |F|^3.$$
From construction of $P^{(N-1)}$ we thus have
$$ |P(\lambda)| \lessapprox |F|^2 \alpha$$
where $P(\lambda) := \lambda \cap P^{(N-1)}$.

Let $L(\lambda)$ denote those lines in $L_1$ which lie in $\lambda$. By construction of $L_1$ we have
$$ |\{ (p,l) \in P(\lambda) \times L(\lambda): p \in l \}| \gtrapprox |F| |L(\lambda)|.$$
On the other hand, from Proposition \ref{wip} we have
$$ |\{ (p,l) \in P(\lambda) \times L(\lambda): p \in l \}| \lessapprox |P(\lambda)|^{1/2} |L(\lambda)|^{3/4} |F|^{1/4} + |P(\lambda)| + |L(\lambda)|.$$
Combining all three estimates and using \eqref{a-small} we obtain
$$ |L(\lambda)| \lessapprox \alpha^2 |F|$$
which is \eqref{no-3space}.

Note in fact that the above argument gives \eqref{no-3space} if $L_1$ is replaced by $L^{(k)}$ for any $1 \leq k \leq N-1$.

We now prove \eqref{no-2plane}.  Following \cite{wolff:xray}, \cite{laba:xray}, we define the \emph{plate number}
$\p_{k}$ for $0 \leq k \leq N-1$ to be the quantity
$$ \p_{k} := \sup_{\pi \in Gr(F^4, 2)} |\{ l \in L_k: l \subset \pi \} |.$$
We observe the bounds
\be{pk-triv}
1 \leq \p_k \lesssim |F|;
\end{equation}
the former bound comes since $L_k$ is non-empty, while the latter bound comes since a 2-plane can contain at most $O(|F|)$ lines in different directions.

Clearly the plate numbers are non-increasing in $k$.  From this, \eqref{pk-triv}, the pigeonhole principle and the fact that $|F|^{1/N} \approx 1$, we can thus find $2 \leq k \leq N-1$ such that
\be{2-pk-const}
\p_{k-1} \approx \p_{k}.
\end{equation}

Fix this $k$.  We can thus find a 2-plane $\pi \in Gr(F^4, 2)$ such that the set
$$ L_k(\pi) := \{ l \in L_k: l \subset \pi \}$$
has cardinality $\p_k$.  

Fix $\pi$, and let $P_k(\pi)$ denote the set 
$$ P_k(\pi) := P_k \cap \pi.$$
By construction of $L_k$, every line in $L_k$ contains $\gtrapprox |F|$ points in $P_k$, thus every line in $L_k(\pi)$ contains $\gtrapprox |F|$ points in $P_k(\pi)$.  In particular we see that
$$ | \{ (p,l) \in P_k(\pi) \times L_k(\pi): p \in l \} | \gtrapprox |F| \p_k.$$
Applying Corollary \ref{easy-cor} we conclude that
$$ |F| \p_k \lessapprox |P_k(\pi)|^{1/2} \p_k + |P_k(\pi)|;$$
from this and \eqref{pk-triv} we thus have
\be{2p-bound}
|P_k(\pi)| \gtrapprox |F| \p_k.
\end{equation}
Let $P'_k(\pi)$ denote the space of all points $p$ in $P_k(\pi)$ such that at least half of all the lines in $\{ l \in L_{k-1}: p \in l \}$ are contained in $L_{k-1}(\pi)$.  We have two cases.

{\bf Case 1 (parallel case): $|P'_k(\pi)| \geq \frac{1}{2} |P_k(\pi)|$.}

In this case we have
\bas
|\{ (p,l) \in P_k(\pi) \times L_{k-1}(\pi): p \in l \}| 
&\geq |\{ (p,l) \in P'_k(\pi) \times L_{k-1}(\pi): p \in l \}| \\
&\geq \frac{1}{2} |\{ (p,l) \in P'_k(\pi) \times L_{k-1}: p \in l \}| \\
&\gtrapprox |P'_k(\pi)| \alpha^{-1} |F|\\
&\gtrapprox |P_k(\pi)| \alpha^{-1} |F|,
\end{align*}
while by definition of $\p_{k-1}$ we have
$$|L_{k-1}(\pi)| \leq \p_{k-1}.$$
Applying Corollary \ref{easy-cor} we thus see that
$$ |P_{k}(\pi)| \alpha^{-1} |F| \lessapprox |P_{k}(\pi)|^{1/2} \p_{k-1} + |P_{k}(\pi)|,$$
which by \eqref{a-small} implies that
$$ |P_{k}(\pi)| |F|^2 \lessapprox \alpha^2 \p_{k-1}^2.$$
But combining this with \eqref{2-pk-const}, \eqref{2p-bound} we obtain 
$$ \p_{k} \gtrapprox \alpha^{-2} |F|^3.$$
But this contradicts \eqref{pk-triv} by \eqref{a-small}.  Hence this case cannot occur.

{\bf Case 2 (transverse case): $|P'_{k}(\pi)| \leq \frac{1}{2} |P_{k}(\pi)|$.}  In this case we have (by a computation similar to Case 1)
$$  |\{ (p,l) \in P_{k}(\pi) \times L_{k}: p \in l; l \not \subset \pi \}| \gtrapprox \alpha^{-1} |P_{k}(\pi)| |F|.$$
Thus, if $L^*_{k-1}$ denotes the lines $l \in L_{k-1}$ which are incident to a point in $P_{k}(\pi)$ but are not contained in $\pi$, then we have
\be{lk-bound}
|L^*_{k-1}| \gtrapprox \alpha^{-1} |P_{k}(\pi)| |F| \gtrapprox \alpha^{-1} |F|^2 \p_{k}
\end{equation}
by \eqref{2p-bound}.

We now use Wolff's hairbrush argument \cite{wolff:kakeya}, \cite{wolff:survey},
as modified to deal with plates in \cite{wolff:xray}, \cite{laba:xray}.
We can foliate $L^*_{k-1}$ as the disjoint union of
$$ L^*_{k-1}(\lambda) := \{ l \in L^*_{k-1}: l \in \lambda \}$$
where $\lambda$ ranges over the 3-spaces containing $\pi$.  For each such $\lambda$, observe from the analogue of \eqref{no-3space} for $L^{(k-1)}$ that 
\be{l-bound}
|L^*_{k-1}(\lambda)| \lessapprox \alpha^2 |F|.
\end{equation}
Also, if we define
$$ P^*_{k-1}(\lambda) := \{ p \in P_{k-1}: p \in \lambda \backslash \pi \}$$
then by construction of $L_{k-1}$, we have
$$ | \{ (p, l) \in P^*_{k-1}(\lambda) \times L^*_{k-1}(\lambda):  p \in l \}|
\gtrapprox |L^*_{k-1}(\lambda)| |F|.$$
Applying Corollary \ref{easy-cor} we obtain
$$ |L^*_{k-1}(\lambda)| |F| \lessapprox |P^*_{k-1}(\lambda)|^{1/2} |L^*_{k-1}(\lambda)| + |P^*_{k-1}(\lambda)|,$$
which by \eqref{l-bound}, \eqref{a-small} implies that
$$ |P^*_{k-1}(\lambda)| \gtrapprox \alpha^{-2} |L^*_{k-1}(\lambda)| |F|.$$
Summing in $\lambda$, we obtain
$$ |P_{k-1}| \gtrapprox \alpha^{-2} |L^*_{k-1}| |F| \gtrapprox \alpha^{-3} |F|^3 \p_k.$$
Since $|P_{k-1}| \lessapprox \alpha |F|^3$ by construction, we obtain $\p_k \lessapprox \alpha^4$, and the claim follows.
\end{proof}

\section{Construction of reguli}\label{np-sec}

We now continue the proof of Theorem \ref{main}.  We begin by refining $P_1$ and $L_1$ a little further. By \eqref{big-incidence} we have
$$ \sum_{l \in L_1} |l \cap P_1| \approx |F|^4.$$
Thus if we set 
$$ L_2 := \{ l \in L_1: |l \cap P_1| \approx |F| \}$$
then by the popularity argument 
$$ \sum_{l \in L_2} |l \cap P_1| \approx |F|^4$$
or equivalently
$$ \sum_{p \in P_1} |\{ l \in L_2: p \in l \}| \approx |F|^4.$$
Thus if we set
$$ P_2 := \{ p \in P_1: |\{ l \in L_2: p \in l \}| \approx \alpha^{-1} |F| \}$$
then by \eqref{p1-bound}, \eqref{high-mult}, and the popularity argument we have
\be{p-bash}
\sum_{p \in P_2} |\{ l \in L_2: p \in l \}| \approx |F|^4.
\end{equation}
In particular we have 
\be{p3-bound}
|P_2| \approx \alpha |F|^3.
\end{equation}

The next task is to generate a large number of frames, and a large number of lines in $L$ incident to the reguli generated by these frames.  As a frame is a fairly complicated combinatorial object (consisting of three lines and a 3-space), we will first begin by counting some simpler objects which eventually will be combined together to form frames.

By \eqref{p-bash} we have
$$ |\{ (p, l) \in P_2 \times L_2: p \in l \}| \approx |F|^4.$$
Since $|L_2| \lessapprox |F|^3$, we thus see from Lemma \ref{cz} that
$$ |\{ (p_1,p_2,l) \in P_2 \times P_2 \times L_2: p_1, p_2 \in l; p_1 \neq p_2 \}| \gtrapprox |F|^8 / |L_2| \approx |F|^5.$$
By definition of $P_2$, we see that for each $(p_1,p_2,l)$ as above, there are
$\gtrapprox \alpha^{-1} |F|$ lines $l_1 \in L_2$ which contain $p_1$ but are distinct from $l$, and similarly there are $\gtrapprox \alpha^{-1} |F|$ lines $l_2 \in L_2$ which contain $p_2$ but are distinct from $l$.  We thus have
$$ |H_0| \gtrapprox \alpha^{-2} |F|^7$$
where
$$ H_0 := \{ (p_1,p_2,l,l_1, l_2) \in P_2 \times P_2 \times L_2 \times L_2 \times L_2: p_1, p_2 \in l; p_1 \neq p_2; p_1 \in l_1; p_2 \in l_2; l \neq l_1, l_2 \}$$
is the space of ``H''-shaped objects.

Let $H_1 \subseteq H_0$ be the set of elements $(p_1,p_2,l,l_1,l_2)$ in $H_0$ such that $l_1$ and $l_2$ are skew.  We claim that
$$ |H_0 \backslash H_1| \lessapprox \alpha^{-1} |F|^6.$$
Indeed, to choose an element $(p_1,p_2,l,l_1,l_2)$ in $H_0 \backslash H_1$ (which is a degenerate ``H'', i.e. a triangle), we first choose $p_1 \in P_2$ (of which there are $\lessapprox \alpha |F|^3$ choices), and then choose the distinct lines $l$, $l_1$ incident to $p_1$ (of which there are $\lessapprox (\alpha^{-1} |F|)^2$ choices).  Since $l_2$ must lie in the 2-plane generated by $l$ and $l_1$, and the lines of $L_1$ point in different directions, there are only $O(|F|)$ choices for $l_2$.  Since $p_2$ is uniquely determined as $p_2 = l \cap l_2$, the claim follows.

From the above bounds and \eqref{a-small} we see that
\be{h1-lower}
|H_1| \gtrapprox \alpha^{-2} |F|^7.
\end{equation}
By construction, if $h = (p_1, p_2, l, l_1, l_2) \in H_1$, then $l_1$ and $l_2$ are skew.  Thus $l_1$ and $l_2$ lie in a unique 3-space $\lambda(h)$, which then must also contain $p_1$, $p_2$, $l$.  

Let $S_0 \subset L_2 \times L_2$ denote the pairs $(l_1, l_2)$ of skew lines in $L_2$.  For each pair $(l_1, l_2) \in S_0$, we define the \emph{connecting set} $C(l_1, l_2) \subset L_2$ to be the set of all lines $l \in L_2$ which are distinct from $l_1$, $l_2$, but intersect both $l_1$, $l_2$ in points $p_1 \in P_2$ and $p_2 \in P_2$ respectively. Observe the identity 
$$ \sum_{(l_1,l_2) \in S_0} |C(l_1,l_2)| = |H_1|.$$
Since $|S_0| \leq |L_2|^2 \lesssim |F|^6$, we thus see from \eqref{h1-lower} that if we define
$$ S_1 := \{ (l_1,l_2) \in S_0: |C(l_1,l_2)| \gtrapprox \alpha^{-2} |F| \},$$
then by the popularity argument
\be{s-bound}
\sum_{(l_1,l_2) \in S_1} |C(l_1,l_2)| \gtrapprox \alpha^{-2} |F|^7.
\end{equation}

If $(l_1,l_2) \in S_1$, we define the set $C^{(3)}(l_1,l_2) \subseteq C(l_1,l_2)^3$ to be the space of all triplets $(l^1, l^2, l^3) \in C(l_1,l_2)^3$ such that the six points $l^i \cap l_j$ for $i=1,2,3$, $j=1,2$ are all disjoint.  

We now use the non-degeneracy property \eqref{no-2plane} to obtain a lower bound for the size of $C^{(3)}(l_1,l_2)$.

\begin{lemma}[Many triple connections between skew lines]\label{triplet}  For any $(l_1, l_2) \in S_1$, we have $|C^{(3)}(l_1,l_2)| \gtrapprox \alpha^{-4} |F|^2 |C(l_1,l_2)|$.
\end{lemma}

\begin{proof}
Fix $l_1$, $l_2$.  We choose $l^1 \in C(l_1, l_2)$ arbitrarily; of course, there are $|C(l_1, l_2)|$ choices for $l^1$.  

Fix $l^1$.  From \eqref{no-2plane} we have
$$ |\{ l^2 \in C(l_1, l_2): l^2 \cap l_1 = l^1 \cap l_1 \}| \lessapprox \alpha^4$$
(since such lines lie in the 2-plane spanned by $l^1 \cap l_1$ and $l_2$.  Similarly if the roles of $l_1$ and $l_2$ are interchanged.  Since $|C(l_1,l_2)| \gtrapprox \alpha^{-2} |F|$, we thus see from \eqref{a-small} that there are $\gtrapprox \alpha^{-2} |F|$ choices for $l^2$ such that $l^2 \cap l_j \neq l^1 \cap l_j$ for $j=1,2$.  

Fix $l^2$.  Arguing as above we see that there are $\gtrapprox \alpha^{-2} |F|$ choices for $l^3$ such that $l^3 \cap l_j \neq l^i \cap l_j$ for $i=1,2$ and $j=1,2$.  The claim follows. 
\end{proof}

From this lemma and \eqref{s-bound} we see that
$$ \sum_{(l_1,l_2) \in S_1} |C^{(3)}(l_1,l_2)| \gtrapprox \alpha^{-6} |F|^9.$$
Observe that if $(l_1, l_2) \in S_1$ and $(l^1,l^2,l^3) \in C^{(3)}(l_1,l_2)$, the various incidence assumptions in the definition of $S_1$ and $C^{(3)}(l_1,l_2)$ force $f := (l^1, l^2, l^3, \lambda)$ to be a frame, where $\lambda$ is the unique 3-space spanned by $l^1$ and $l^2$.  Observe that $l_1, l_2$ both lie in $L_2 \cap L(f)$.  Thus, if $\FRAME_0$ denotes the space of all frames generated in this manner, then
\be{f-big}
\sum_{f \in \FRAME_0} | L_2 \cap L(f) |^2 \gtrapprox \alpha^{-6} |F|^9.
\end{equation}
Let $f \in \FRAME_0$.  Since the lines in $L(f)$ are contained in a regulus, they have finite overlap.  Since each line in $L_2$ contains $\approx |F|$ points in $P_1$ by construction, we thus see that\footnote{One could also obtain this bound using Corollary \ref{easy-cor} and the crude bound $|L_2 \cap L(f)| \leq |L(f)| \lesssim |F|$.}
$$ |P_1 \cap r(f)| \gtrsim |F| |L_2 \cap L(f)|$$
so by \eqref{f-big} we have
$$ \sum_{f \in \FRAME_0} |P_1 \cap r(f)|^2 \gtrapprox \alpha^{-6} |F|^{11}.$$
By \eqref{high-mult}, each point in $P_1$ is incident to $\approx \alpha^{-1} |F|$ lines in $L_1$.  Thus we have
$$ \sum_{f \in \FRAME_0} |\{ l \in L_1: l \cap r(f) \cap P_1 \neq \emptyset \}|^2 \gtrapprox \alpha^{-8} |F|^{13}.$$
We observe the cardinality bound
\be{f0}
|\FRAME_0| \lessapprox \alpha^2 |F|^7.
\end{equation}
Indeed, to choose a frame $(l^1,l^2,l^3,\lambda)$ in $\FRAME_0$, we observe that there are $O(|L_1|^2) = O(|F|^6)$ choices for the skew pair $(l^1,l^2)$.  This determines $\lambda$, and then by \eqref{no-3space} we thus see that there are $O(\alpha^2 |F|)$ choices for $l^3$, and \eqref{f0} follows.  In particular, if we define
\be{f1-def}
\FRAME_1 := \{ f \in \FRAME_0: |\{ l \in L_1: l \cap r(f) \cap P_1 \neq \emptyset \}| \gtrapprox \alpha^{-5} |F|^3 \}
\end{equation}
then by the popularity argument 
\be{frame-1}
\sum_{f \in \FRAME_1} |\{ l \in L_1: l \cap r(f) \cap P_1 \neq \emptyset \}|^2 \gtrapprox \alpha^{-8} |F|^{13}.
\end{equation}
Since the summand on the left-hand side can be crudely bounded by $|L_1|^2 = O(|F|^6)$, we thus have the crude bound
\footnote{The bounds on $|\FRAME_1|$, and later on $|\FRAME_2|$, $|\FRAME_3|$, might not be best possible, however an improvement on this part of the argument does not directly improve the gain $\gain$.}
\be{f1-lower}
|\FRAME_1| \gtrapprox \alpha^{-8} |F|^7
\end{equation}
(compare with \eqref{f0}).

For any frame $f \in \FRAME_1$, there are only $O(|F|^3)$ possible orientations for $\lambda(f)$.  By \eqref{f1-lower} and the pigeonhole principle, there therefore exists a 3-space $\lambda_0 \in Gr(F^4, 3)$ such that
\be{r-large}
|\FRAME_2| \gtrapprox \alpha^{-8} |F|^4
\end{equation}
where
$$ \FRAME_2 := \{ f \in \FRAME_1: \lambda(f) \hbox{ is a translate of } \lambda_0 \}.$$
Fix this $\lambda_0$.  Let $\FRAME_3$ be a maximal subset of $\FRAME_2$ such that the reguli $\{ r(f): f \in \FRAME_3 \}$ are all distinct.  Since each $r(f)$ contains at most $O(|F|)$ lines, each regulus can arise from at most $O(|F|^3)$ frames.  We thus see from \eqref{r-large} that
\be{r-large-2}
|\FRAME_3| \gtrapprox \alpha^{-8} |F|.
\end{equation}
From \eqref{f1-def} we have 
$$ \sum_{f \in \FRAME_3} | \{ l \in L_1: l \cap r(f) \cap P_1 \neq \emptyset \} | \gtrapprox \alpha^{-5} |F|^3 |\FRAME_3|.$$
From \eqref{a-small}, \eqref{r-large-2} the right-hand side is $\gg |F|^3 \gtrsim |L_1|$.  Thus we can use Lemma \ref{cz}, and obtain
$$ \sum_{f_1,f_2,f_3 \in \FRAME_3: f_1,f_2,f_3 \hbox{ distinct}} | \{ l \in L_1:
l \cap r(f_i) \cap P_1 \neq \emptyset \hbox{ for } i=1,2,3 \}| \gtrapprox \alpha^{-15} |F|^3 |\FRAME_3|^3.$$
From the pigeonhole principle, we may thus find distinct frames $f_1, f_2, f_3$ in $\FRAME_3$ such that
\be{l1}
|L_*| \gtrapprox \alpha^{-15} |F|^3,
\end{equation}
where $L_* \subseteq L_1$ is the collection of lines
$$ L_* := \{ l \in L_1: l \cap r(f_i) \cap P_1 \neq \emptyset \hbox{ for } i=1,2,3 \}.$$

Now we consider the problem of obtaining upper bounds on $|L_*|$.  The crude upper bound of $|L_1| \sim |F|^3$ is clearly not enough to obtain a contradiction.  However, thanks to the three regulus lemma we can improve this bound by about $|F|$:

\begin{proposition}\label{up-bound}  We have
\be{up-l}
|L_*| \lessapprox |F|^{2+\gain}.
\end{equation}
\end{proposition}

\begin{proof}
If the 3-spaces $\lambda(f_1)$, $\lambda(f_2)$, $\lambda(f_3)$ are disjoint, then this follows directly from Corollary \ref{3-reg} and \eqref{p0-small}.

By symmetry it remains to consider the case when $\lambda(f_1)$ and $\lambda(f_2)$ (for instance) are equal.  Then the lines in $L_*$ must either be parallel to $\lambda(f_1)$, or else intersect $P_1 \cap r(f_1) \cap r(f_2)$.  There are at most $|F|^2$ lines in the first category (in fact there are far fewer, thanks to \eqref{no-3space}).  In the second category, we observe that $r(f_1) \cap r(f_2)$ is at most one-dimensional (since $r(f_1)$, $r(f_2)$ are irreducible and distinct) and hence has cardinality $O(|F|)$.  On the other hand, by \eqref{high-mult} every point in $r(f_1) \cap r(f_2) \cap P_1$ is incident to $\approx \alpha^{-1} |F|$ lines in $L_1$.  Thus we certainly have $\lessapprox |F|^{2+\gain}$ incidences in this case as well. 
\end{proof}

Combining \eqref{up-l} with \eqref{l1} we obtain 
$$ \alpha \gtrapprox |F|^\gain$$
and hence by \eqref{p3-bound}
$$ |P_0| \gtrapprox |P_2| \approx \alpha |F|^3 \gtrapprox |F|^{3 + \gain}$$
as desired. 
This concludes the proof of Theorem \ref{main}.

\section{Remarks}\label{remarks-sec}

It seems likely that this Theorem can be generalized in several ways.  The exponent $\gain$ is probably not sharp, and also the result should have extension to other dimensions\footnote{In dimensions 5 and higher there are other, more ``arithmetic'' arguments which give slight improvements to $|F|^{(n+2)/2}$ for Besicovitch sets; see \cite{borg:high-dim}, \cite{KT}, \cite{gerd:kakeya}, \cite{rogers}, \cite{kt:kakeya2}, \cite{kt-survey}.  Nevertheless, if one can make an improvement of the order of $\gain$ in, say, five dimensions by these ``geometric'' techniques, this will be quite competitive with the results in, say, \cite{kt:kakeya2}.  In the Euclidean setting one can improve the bound $(n+2)/2$ in all dimensions $n \geq 3$ by a small number ($10^{-10}$) for the \emph{upper Minkowski dimension} problem for Besicovitch sets (\cite{katzlaba}, \cite{laba:xray}, \cite{laba:medium}), but this argument seems special to the upper Minkowski problem and does not directly impact the finite field question.} although the amount of algebraic geometry needed to reproduce the above argument in higher dimensions seems non-trivial.  Also, the argument can probably be extended to obtain an estimate on the Kakeya maximal function for finite fields (see \cite{gerd:kakeya}).  In principle, the finite field results should also extend to the Euclidean setting $\R^n$, but there are some unpleasant technical difficulties\footnote{In the finite field case one is aided considerably by the fact that dimensions must be integer; for instance, the intersection of two lines is either empty, 0-dimensional (a point), or 1-dimensional (a line).  In the ($\delta$-discretized) Euclidean case there is a continuum of cases: two distinct $1 \times \delta$ tubes can intersect in a set of length $\sim \delta/\theta$, where $\delta < \theta < 1$ is the angle between the two tubes.  This introduces a new dyadic parameter $\theta$ into the analysis (measuring the degeneracy of the angle), and often the cases of small $\theta$ and large $\theta$ need to be treated separately.  See e.g. \cite{wolff:kakeya}.   Here we have more complicated algebraic objects, such as the variety \eqref{W-def}, and to capture the possible degeneracies of this object seems to require a large number of additional dyadic parameters.  It is possible that various rescaling arguments, such as the two-ends and bilinear reductions mentioned above, may be used to reduce the number of such parameters, but the extension of this argument to the Euclidean case still appears to be quite non-trivial.} in doing so due to the presence of near-degenerate reguli, etc. in $\R^n$.  (These type of difficulties cause considerable complication in such papers as \cite{schlag:kakeya},
although some of this difficulty could perhaps be alleviated using the ``two-ends'' reduction in \cite{wolff:kakeya} and the ``bilinear reduction'' in \cite{tvv}.).  We will not pursue these matters.

\end{document}